\documentclass[11pt]{amsart}

\newif\ifpdf
    \ifx\pdfoutput\undefined
    \pdffalse 
    \else
    \pdfoutput=1 
    \pdftrue
    \fi

    \ifpdf
    \usepackage[pdftex]{graphicx}
    \else
    \usepackage{graphicx}
    \fi
\usepackage{amsmath,amssymb,amsthm,epsfig}

\newcommand\G{\Gamma}
\newcommand\p{p_1,p_2, \cdots ,p_r}
\newcommand\pexp{p_1^{e_1},p_2^{e_2}, \cdots ,p_r^{e_r}}
\newcommand\g{\gamma}
\newcommand\A{{\mathbb H}^2}
\newcommand\R{\mathbb{R}}
\newcommand\Q{\mathbb{Q}}

\newcommand\Z{\mathbb{Z}}

\newcommand\Qp{\mathbb{Q}_p}

\newcommand\F{\Z [\frac{1}{p}]}
\newcommand\PZp{PSL_2( \F)}

\newcommand\PZn{PSL_{2}(\Z[\frac{1}{n}])}
\newcommand\PZm{PSL_{2}(\Z[\frac{1}{m}])}

\newcommand\PQp{PSL_2(\Qp)}
\newcommand\PR{PSL_2(\R)}
\newcommand\semi{\rtimes}
\newcommand\s{(n_1,n_2, \cdots n_k)}
\newcommand\sprime{(m_1,m_2, \cdots ,m_l)}
\newcommand\pt{\Pi_{i=1}^k T_i}
\newcommand\pti{\Pi_{i=1}^k T^i}
\newcommand\point{(t_1,t_2, \cdots ,t_k)}
\newcommand\ptl{\Pi_{i=1}^l T_i}
\newcommand\ptr{\Pi_{i=1}^r T_i}
\newcommand\ptp{\Pi_{i=1}^k T'_i}

\newtheorem{theorem}{Theorem}[section]

\newtheorem{lemma}[theorem]{Lemma}
\newtheorem{corollary}[theorem]{Corollary}

\input{epsf.tex}

\begin{document}

\title{The Large Scale Geometry of Some Metabelian Groups}
\author{Jennifer Taback and Kevin Whyte}
\thanks{The first author acknowledges partial support from the NSF.  The second author acknowledges partial support from a University at Albany/UUP Individual Development Award.}

\begin{abstract}  We study the large scale geometry of the upper triangular
subgroup of $\PZn$, which arises naturally in a geometric context.  We
prove a quasi-isometry classification theorem and show that these
groups are quasi-isometrically rigid with infinite dimensional
quasi-isometry group.  We generalize our results to a larger class
of groups which are metabelian and are higher dimensional analogues of the
solvable Baumslag-Solitar groups $BS(1,n)$.
\end{abstract}

\maketitle
\section{Introduction}

We consider quasi-isometries of the upper triangular subgroup
$\G_n$ of $\PZn$.  These groups arise in a geometric way because
they are subgroups of both $\PR$ and  $\PQp$, for all $p$ dividing $n$.
Both $\PR$ and $\PQp$ act on their respective Bruhat-Tits buildings;
for $\PQp$ this building is
a regular $p+1$ valent tree, and for $\PR$ it is $\A$.
Then $G=\PR \times \Pi_{p_i | n} PSL_2(\Q_{p_i})$
has an induced action on $\A \times \pt$, where $T_i$ is the Bruhat-Tits
building of $PSL_2(\Q_{p_i})$.  This action is properly
discontinuous and has cofinite volume, but its restriction to $\G_n$ is no longer cofinite. However, the induced
action of $\G_n$ on the product of trees is cocompact; the quotient
is a $k$-torus. The stabilizer of any point is an infinite
cyclic group which acts parabolically on $\A$. Thus $\G_n$ has a
decomposition as a $k$ dimensional complex of groups \cite{BH}.

The upper triangular subgroup $\G_n$ arises naturally as the stabilizer
of a point at infinity under the action of $G$
on $\A \times \pt$.  For $n$ prime, this group of
upper triangular matrices is isomorphic to the solvable
Baumslag-Solitar group $BS(1,p^2) = \langle a,b | aba^{-1} = b^{p^2} \rangle$,
and our results on
quasi-isometries and rigidity generalize the results of \cite{FM}.
In this case, the rigidity of the groups $\G_n$ should be useful for
understanding the groups $PSL_2(\Z[\frac{1}{n}])$, analogously to how the
results of Farb and Mosher are used in \cite{T}.

The upper triangular groups $\G_n$ are also basic examples of metabelian groups,
fitting into the short exact sequence
$$ 1 \to \Z[\frac{1}{n}] \to \G \to \Z^k \to 1.$$

In the sections below, we describe geometric
models for these groups as warped products of $\R$ with the product of trees on which $\G_n$ acts.  This
identifies $\G_n$ as a cocompact lattice in the isometry group
$\R \semi ( {Sim}(\Q_{m_1}) \times \cdots
\times {Sim}(\Q_{m_k}))$
of
this model space, where ${Sim}(\Q_m)$ is the group of
similarities of the $m$-adic rationals.  We also describe the  group of all
self quasi-isometries of $\G_n$ and classify them up to quasi-isometry.

Our results rely on the technology available for groups
acting on trees.  However, products of trees are substantially
more complicated than trees.  For example, a group which acts
freely on a tree is free, while groups which act freely on a
product of trees need not be products of free groups.  Such groups
can, in fact, be simple \cite{BM}.

Our results generalize immediately to a larger class of groups which do
not arise as nicely in a geometric context but are interesting nonetheless.
This larger class of groups generalizes the
solvable Baumslag-Solitar groups $BS(1,n) = \langle  a,b | aba^{-1} =
b^n \rangle$. Let $S = \s$ where $(n_i,n_j) = 1$ when $i \neq j$, and
define $\G = \G(S)$ by
$$\G = \G(S) = \langle a_{1},\ldots,a_{k},b | a_{i}^{-1}ba_{i}=b^{n_{i}}, \ a_ia_j = a_ja_i, \ i \neq j \rangle.$$
These groups are $k+1$
dimensional  metabelian groups, fitting into a short exact sequence
$$ 1 \to A \to \G \to \Z^{k} \to 1$$
where the map onto $\Z^{k}$ is given by sending the $\{a_{i}\}$ to a
basis and sending $b$ to $0$.  The kernel, $A$, is normally generated
by $b$ and is an infinitely generated abelian group.   Thus these
groups provide natural examples of finite type solvable groups which
are not polycyclic.

The groups $\G_n$ are also of the above form.
Namely, let $n = p_1^{e_1}p_2^{e_2} \cdots p_k^{e_k}$ where the $p_i$ are
distinct primes.  Then $\G_n$ is isomorphic to
$\G(p_{1}^{2e_{1}},\ldots, p_{k}^{2e_{k}})$, where the isomorphism is
given by:

$$a_{i} \mapsto
\left(\begin{array}{ccc}
 p_{i}^{e_{i}} & 0 \\
0 & {p_i}^{-e_i} \\
\end{array}\right) $$

$$b \mapsto
\left(\begin{array}{ccc}
 1 & 1 \\
0 & 1 \\
\end{array}\right). $$

The decomposition of $\G_n$ into a $k$-dimensional complex of groups
can be generalized to the groups $\G(S)$.  Indeed,
the presentation given is that of a $k$-torus of infinite cyclic
groups, generalizing the fact that all the Baumslag-Solitar groups
are HNN extensions of $\Z$.  This decomposition is fundamental to
our study of the geometry of these groups.  The groups $\G(S)$
have geometric models analogous to those of the $\G_n$.  As a result,
our quasi-isometry classification and rigidity results immediately
generalize to this larger class of groups.  We are able to identify $\G(S)$
as a cocompact lattice in the isometry group of
the model space, describe its quasi-isometry group, and classify
these groups up to quasi-isometry.

\subsection{Statement of Results}

Let $\G_n$ be the upper triangular subgroup of $\PZn$ and $X_n$ the model space for
$\G_n$ which is quasi-isometric to $\G_n$ and constructed below in \S
\ref{sec:models}.

\begin{theorem}[Quasi-isometry classification]
\label{thm:classification}
Let $\G_n$ be the upper triangular subgroup of $\PZn$, and $\G_m$ the upper triangular subgroup of $\PZm$.  If $n = p_1^{e_1} \cdots p_k^{e_k}$ and $m = q_1^{f_1} \cdots q_l^{f_l}$ where $\{p_i\}$ and $\{q_j\}$ are sets of distinct primes, then $\G_n$ and $\G_m$ are quasi-isometric iff $k=l$ and for $i=1,2, \cdots ,k$, after possibly reordering, $p_i = q_i$.
\end{theorem}

\begin{theorem}[Quasi-isometry group]
\label{thm:QIgroup}
Let $n=p_1^{e_1} \cdots p_k^{e_k}$ where all $p_i$ are distinct primes.
The quasi-isometry group, $QI(\G_n)$, is isomorphic to the product
$$ Bilip(\R) \times Bilip(\Q_{p_1}) \times \cdots
\times Bilip(\Q_{p_k}).$$
\end{theorem}

\begin{theorem}[Cusp group rigidity]
\label{thm:rigidity}
If $\G'$ is a finitely generated group which is quasi-isometric to $\G_n$ then
there is a finite normal subgroup $F$ of $\G'$ so that $\G'/F$ is commensurable to $\G_n$,
meaning that
$\G'/F$ and $\G_n$ have isomorphic subgroups of finite index.
\end{theorem}

When we replace $\G_n$ by the more general group $\G(S)$ defined above, where all
elements in $S$ are pairwise relatively prime, we obtain the following generalizations
of the above theorems.

\begin{theorem}
\label{thm:classification2} Consider the sets $S_1 = \s$ and $S_2
= \sprime$ with $(n_i,n_j) = (m_i,m_j) = 1$ for $i \neq j$. Define
$\G_1 = \G(S_1)$ and $\G_2 = \G(S_2)$. The groups $\G_1$ and $\G_2$
are quasi-isometric iff $k=l$ and for $i=1,2, \cdots k$, after
possibly reordering, each $n_i$ is a rational power of $m_i$.
\end{theorem}

\begin{theorem}
\label{thm:QIgroup2}
Let $S= \s$.  The quasi-isometry group $QI(\G(S))$ is isomorphic to the
product
$$ Bilip(\R) \times Bilip(\Q_{n_1}) \times \cdots
\times Bilip(\Q_{n_k}).$$
\end{theorem}

\begin{theorem}
\label{thm:rigidity2}
Let $\G'$ be any finitely generated group quasi-isometric to $\G(S)$.
There are integers $m_1, m_2, \cdots, m_k$,
with each $m_i$ a rational power of $n_i$, and a finite normal subgroup $F$ of $\G'$ so that
$\G'/F$ is isomorphic to a cocompact lattice in $Iso(X(m_1,\cdots,m_k))= \R \semi ({Sim}(\Q_{m_1})
\times \cdots \times {Sim}(\Q_{m_k}))$.
\end{theorem}

\subsection{Outline of the Proofs}

The key to all of our results is understanding the self quasi-isometries of the model
space $X=X_n$ for $\G_n$, and in general for $\G(S)$, constructed in \S \ref{sec:models}.
This model space is the warped product of $\R$ and a product
of trees $\pt$.
We begin with a definition crucial to understanding the following
outline of the proofs, and refer the reader to \S \ref{sec:prelim} for
additional definitions.
Throughout, let $f: X \rightarrow X$ be any quasi-isometry.

When considering points in $\pt$ it is important to define a notion of
{\em height} on each tree $T_i$.
Fix a basepoint $(t_1,t_2, \cdots ,t_k) \in \pt$.  The height
of a vertex $t \in T_i$ is the height change between $t$ and the
$i$-th coordinate $t_i$ of the basepoint.  Extend this notion to a
height function $h_i$ on each tree $T_i$ through linear interpolation
along the edges.  The metric on $X$ is then given by a warped product of $\R$ and $\pt$
where on each tree $T_i$ the warping function is given by $e^{-h_i}$.

In the following outline, as in the majority of the paper, we only consider the groups
$\G_n$.  The similarities in the construction of the model spaces for the groups $\G_n$
and $\G(S)$ ensure that the generalizations of the proofs are immediate. Let $n= \pexp$ where the $p_i$ are distinct primes.

\begin{itemize}

\item
{\bf Warped product structure is preserved.}
We first show that any quasi-isometry preserves, up to bounded distance, the
{\em horocycles}, i.e. the subsets of the form $\R \times (t_1,\cdots,t_k)$. In other words,
that there is a quasi-isometry $\bar{f}$ of the product of trees $T_1 \times \cdots
\times T_k$, so that $f$ and $\bar{f}$ commute with the projection $X \to T_1 \times
\cdots \times T_k$. Results of \cite{KL} imply that, up to permuting the
factors, $\bar{f}$ splits as a product of quasi-isometries $\bar{f_i}$ of the trees $T_i$.

\item
{\bf Quasi-isometries are {\em almost height translations} on the tree
  factors.}
The geometry of the space $X$ restricts the quasi-isometries $\bar{f_i}$. The
warping function can be reconstructed as the (logarithm of the) amount of stretching induced by
closest point projection between the horocycles.  This splits as a sum of functions, $h_i$,
on each of the trees.  The quasi-isometries $\bar{f_i}$ preserve these warping functions,
in the sense that $h_i(\bar{f_i}(x))-h_i(\bar{f_i}(y))$ differs from $h_i(x)-h_i(y)$ by
a uniformly bounded amount.  We call such quasi-isometries {\em almost height translations}.
In \cite{FM}, the group of almost height translations of $T_n$ is identified as $Bilip(\Q_n)$.

\item
{\bf $f$ induces a bilipschitz homeomorphism of $\R$.} This shows
that the group of quasi-isometries of $T_1 \times \cdots \times
T_k$ which quasi-preserve the warping function is $Bilip(\Q_{p_1})
\times \cdots \times Bilip(\Q_{p_r})$. All of these
quasi-isometries extend to quasi-isometries of $X$.  The
quasi-isometries of $X$ which induce the identity on $T_1 \times
\cdots \times T_r$ induce a bilipschitz homeomorphism of $\R$.
This allows us to identify the quasi-isometry group of $X$, and
prove theorem \ref{thm:QIgroup} (quasi-isometry group).

\item
{\bf These methods hold for quasi-isometries between $\G_n$ and
  $\G_m$.}
Consider a quasi-isometry $f:\G_n \rightarrow \G_m$. Using the above
methods again shows that $f$ induces a bilipschitz homeomorphism
of $\R$ and a quasi-isometry on the product of trees which is a
bounded distance from a product quasi-isometry. Theorem
\ref{thm:classification} (quasi-isometry classification) now
follows by combining results of \cite{FM} and \cite{W1}.

\item
{\bf Quasi-actions.} Understanding the quasi-isometries of $X$ lets us understand groups quasi-isometric to $\G$
via the quasi-action principle.  Suppose $\G'$ is quasi-isometric to $\G_n$ (and hence to $X$),
and let $f:\G' \to X$ be a quasi-isometry.  For every $\g' \in \G'$ we get a quasi-isometry
of $X$ by $x \mapsto f(\g' f^{-1}(x))$.  These quasi-isometries all have uniform constants, and
compose, up to bounded distance, according to the multiplication table of $\G'$.  In other
words, $\G'$ quasi-acts on $X$, and therefore gives an almost height translation action on each of
the $T_i$.

\item
{\bf Obtaining similarity actions on $\Q_n$.}
According to \cite{MSW}, these almost height translation actions are equivalent, via
a quasi-isometry $T_i \to T'_i$, to a height translation action on trees $T'_i$.  In terms
of the $p_i$-adics, this says that there is some $q_i$ so that the bilipschitz action
of $\G'$ on $\Q_{p_i}$ is bilipschitz equivalent to a similarity action on $\Q_{q_i}$.
Similarly, the bilipschitz action on $\R$ is equivalent to an affine action on $\R$.  Further,
the uniformity of the quasi-isometry constants implies that the expansion factor of the
affine action on $\R$ is the inverse of the product of the factors from the similarity
actions on the $\Q_{q_i}$.  This shows $\G'$ is a lattice in the subgroup of ${Aff}(\R)\times
{Sim}(\Q_{q_1}) \times \cdots \times {Sim}(\Q_{q_k})$ which satisfies this condition.  This
subgroup is $\R \semi {Sim}(\Q_{q_1}) \times \cdots \times {Sim}(\Q_{q_k})$, and can be
identified as the isometry group of a complex $X'$, proving theorem \ref{thm:rigidity}.

\end{itemize}

\section{Preliminaries}
\label{sec:prelim}

\subsection{Definitions and Notation}

We begin with the definition of a quasi-isometry.

\medskip
\noindent
{\bf Definition.}  Let $K \geq 1$ and $C \geq 0$.
A {\it $(K,C)$-quasi-isometry} between metric
spaces $(X, d_X)$ and $(Y,d_Y)$ is a
map $f:X \rightarrow Y$ satisfying:

\medskip
\noindent
1. $\frac{1}{K} d_X(x_1,x_2) - C \leq d_Y(f(x_1),f(x_2)) \leq K
d_X(x_1,x_2) + C$ for all $x_1,x_2 \in X$.

\medskip
\noindent
2. For some constant $C'$, we have $Nbhd_{C'}(f(X)) = Y$.

\bigskip

We will assume that our quasi-isometries have been changed by a bounded amount using
the standard ``connect-the-dots'' procedure to be continuous.
(See, for example, \cite{SW}.)
A quasi-isometry has a {\em coarse inverse}, i.e. a quasi-isometry $g:Y
\rightarrow X$ so that $f \circ g$ and $g \circ f$ are a bounded
distance from the appropriate identity map in the sup norm.
A map satisfying $1.$ but not $2.$ in the definition above is called a
{\em quasi-isometric embedding}.

We define the {\em quasi-isometry group} $QI(X)$ of a space $X$  to
be the collection of all self quasi-isometries of $X$, identifying
those which differ by a bounded amount in the sup norm.

Given a group $G$ and a metric space $X$, a {\em quasi-action} of $G$ on $X$ associates to each $g \in G$ a quasi-isometry of $X$, i.e. $A_g: X \rightarrow X$, subject to certain conditions.  This map is defined by $A_g(x) = g \cdot x$, and the collection of these maps has uniform quasi-isometry constants, so that $A_{Id} = Id_X$ and $d_{sup}(A_g \circ A_h,A_{gh})$ is bounded independently of $g$ and $h$.

\subsection{Previous results}

The following theorems will be referred to repeatedly in \S
\ref{sec:proofs}.  We state them below for easy reference.

\subsubsection{Rigidity of Baumslag-Solitar groups}
\label{sec:FM}

Since the geometry of the group $\Gamma_n$ is so dependent on its
various Baumslag-Solitar subgroups, we will often refer to the
following classification and rigidity results for the solvable
Baumslag-Solitar groups due to Farb and Mosher.

\begin{theorem}[\cite{FM} Theorem 7.1]
\label{thm:BS:1}
For integers $m,n \geq 2$, the groups $BS(1,m)$ and $BS(1,n)$ are
quasi-isometric if and only if they are commensurable.  This happens
if and only if there exist integers $r,j,k>0$ such that $m=r^j$ and $n
= r^k$.
\end{theorem}

\begin{theorem}[\cite{FM} Theorem 8.1]
\label{thm:BS:2}
The quasi-isometry group of $BS(1,n)$ is given by the following isomorphism:
$$QI(BS(1,n)) \cong Bilip(\R) \times Bilip(\Q_n).$$
\end{theorem}

\subsubsection{Products of trees and groups acting on products of
  trees}

A major step in the proofs below is to show that a quasi-isometry
$f:\G_1 \rightarrow \G_2$ induces a map on the product of trees on
which each group acts.  Once this is accomplished, we use the
following result of Kleiner and Leeb to show that our map is
uniformly close to a product of quasi-isometries.

\begin{theorem}[\cite{KL} Theorem 1.1.2]
\label{thm:KL} Let $T_i$ and $T_i'$ be irreducible thick Euclidean
Tits buildings with cocompact affine Weyl group. Let $X = {\mathbb
E}^n \times \Pi_{i=1}^k T_i$ and $X' = {\mathbb E}^{n'} \times
\Pi_{i=1}^{k'} T_i'$ be a metric products.  Then for all $K,C > 0$
there exist $K', \ C', \ D'$ so that the following holds: If $f: X
\rightarrow X'$ is a $(K,C)$-quasi-isometry, then $n=n'$, $k=k'$
and there are $(K',C')$-quasi-isometries $f_i : T_i \rightarrow
T_j$ so that $d(p \circ f, \Pi_{i=1}^k f_i \circ p) \leq D'$ where
$p$ is the projection map.
\end{theorem}

The following result of \cite{MSW} will be needed for the proof of
rigidity of the groups $\G$.  It applies to {\em bushy} trees, meaning that each vertex is a uniformly bounded distance from a vertex having at least three unbounded complementary components.  In addition we require {\em bounded valence}, meaning that vertices have uniformly finite bounded valence.  All of the trees in the discussion below satisfy these properties.

\begin{theorem}[\cite{MSW}]
\label{thm:MSW}
If $G \times T \rightarrow T$ is a quasi-action of a group $G$ on a bounded valence, bushy tree $T$, then there is a bounded valence, bushy tree $T'$, an isometric action $G \times T' \rightarrow T'$, and a quasi-isometry $f: T' \rightarrow T$ which intertwines the actions of $G$ on $T'$ and the quasi-action of $G$ on $T$ to within a uniformly bounded distance.
\end{theorem}

\section{The Geometric Models}
\label{sec:models}

To illustrate the geometry of $\Gamma_n$, and $\G(S)$ in general, we describe a metric
$(k+1)$-complex $X$ quasi-isometric to $\Gamma_n$, i.e. on which
$\Gamma_n$ acts properly discontinuously and cocompactly by
isometries. We begin with the simplest case of the upper triangular
subgroup $\G_n$ of $\PZn$.  We then describe the geometry of the more general
groups $\G(S)$.  For all of these groups, the complex $X$ is a warped
product of $\R$ with a product of trees on which the group acts.
When the $n_i$ are not relatively prime, the group $\G(S)$ does not
act on a product of trees, and we do not consider this case here.

\begin{figure}
\includegraphics[width=2in]{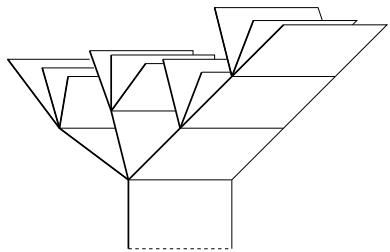}\\
\caption{The geometric model of the solvable Baumslag-Solitar group $BS(1,3)$, which is topologically a warped product of a tree and $\R$. \label{fig:bs}}
\end{figure}

First recall that the Baumslag-Solitar groups $BS(1,n) = \langle
a,b | aba^{-1} = b^n \rangle$, for integral $n \geq 2$, for integral $n \geq 2$, act properly discontinuously and
cocompactly by isometries on a metric $2$-complex we denote $Y_n$.
This complex $Y_n$ is topologically the product $T \times \R$
where each vertex of the tree has $1$ incoming edge and $n$
outgoing edges. Metrically we define a height function on $T$ so
that if $l \subset T$ is a line on which the height function is
strictly increasing, then $l \times \R$ is metrically a hyperbolic plane (\ref{fig:bs}).
See \cite {FM} for a more detailed construction of this complex, and figure .

\subsection{The geometric model of $\G_n$}

We give the most comprehensive description of the model space $X$
in this case because the trees on which $\G_n$ acts are easier
to understand than the trees on which $\G(S)$ acts.
We present several ways to understand the complex
$X$.

When $p$ is prime, the group $BS(1,p)$, acts on the Bruhat-Tits tree $T_p$ associated to
$PSL_2(\Q_p)$. This is not true for $BS(1,n)$ when $n$ is not prime. We
will describe the $BS(1,n)$ tree below.

\begin{figure}
\includegraphics[width=1in]{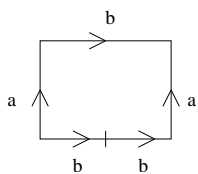}\\
\caption{The ``horobrick" building block for the geometric model of $BS(1,2)$. \label{fig:brick1}}
\end{figure}

Assume $p$ is prime, and consider the geometric model $Y_p$ of
$BS(1,p)$. Let $(x,y)$ be the coordinates on the upper half space
model of hyperbolic space, where $y > 0$. One can also view $Y_p$
as built from the ``horobrick'' with $0 \leq x \leq n$ and $1 \leq
y \leq p$. The vertical sides of this brick have length $\log p$.
In the Cayley graph of $BS(1,p)$ this horobrick has the form
given in figure \ref{fig:brick1}.

If $n = \pexp$, then $\G_n$ acts on the product
of the trees $\pt$, where $T_i$ is the tree on which $BS(1,p_i)$
acts, i.e. it has $1$ incoming edge at each vertex and $p_i$
outgoing edges.  The complex $X$ is the same warped product of
$\pt$ with $\R$ as we saw above for $BS(1,p)$.

\begin{figure}
\includegraphics[width=1.5in]{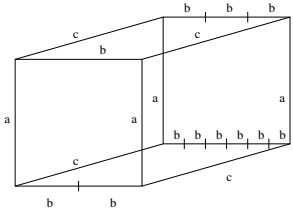}\\
\caption{The analogous building block for $\Gamma(2,3)$. \label{fig:brick2}}
\end{figure}

Analogously for $\Gamma_n$, there is an $(k+1)$-dimensional building block used
to construct the complex $X$, whose $1$-skeleton is the Cayley
graph of $\Gamma$.
When $n$ is a product of two primes, an examples of this block is given in figure \ref{fig:brick2}.
It is not difficult to see that the correct branching occurs when
these blocks are arranged so as to form the appropriate
Baumslag-Solitar subcomplexes.  In general the $(k+1)$-dimensional building block will be an
$(k+1)$-cube, with appropriate edge labels in terms of the
generators of $\Gamma(S)$.  We refer to the horocycle along which the
sheets meet as {\em branching horocycles}.

A second way of understanding the complex $X$ is in terms of
some of its special subspaces.
Let $n = p_1^{e_1}p_2^{e_2} \cdots p_k^{e_k}$ where the $p_i$ are
prime.  Then $$\G_n \cong \langle a_{1},\ldots,a_{k},b | a_{i}^{-1}ba_{i}=b^{p_i^{2e_i}}, \ a_ia_j = a_ja_i, \ i \neq j \rangle.$$
We consider in particular two types of subspaces of $X$:

\begin{itemize}
\item
$Y_{p_i}$, corresponding to $BS(1,p_i^{2e_i})$ generated by $a_i$ and $b$
in the presentation above

\medskip
\item
$\Z^l$, for $1 \leq l \leq k$, generated by $l$ distinct generators
$a_i$ in the presentation above.

\end{itemize}

Notice that the $BS(1,p_i)$ subgroups of $\Gamma$ all share the
generator $b$; In $X$ this means that the subcomplexes $Y_{p_i}$
for $i = 1,2, \cdots ,k$, are joined along branching horocycles.
Namely, consider a subcomplex $Y_{p_i}$ of $X$. At each branching
horocycle of $Y_{p_i}$ there is a copy of $Y_{p_j}$ for all $j
\neq i$ attached along that horocycle. The same is true for every
branching horocycle of those $Y_{p_j}$ and the process continues.

For any point $x \in X$, there is a $Y_{p_i}$ subspace for each
$i = 1,2, \cdots ,n$ in $X$ which contains $x$. For each $i$, the
set $\{ a_i^m \cdot x | m \in \Z \}$ is a line in the Cayley graph
of $\Gamma$ which is the $1$-skeleton of $X$. These lines form the
axes of a $\Z^n$ subspace of $X$. The orbit of $x$ under the group
generated by the entire collection $\{a_i\}$ is a $\Z^n$ subspace
of $X$, for any $x \in X$. The $\Z^l$ subspaces for $l < k$ are
contained in the $\Z^n$ subspaces and are formed by taking the
orbit of $x \in X$ under the group generated by a subset of $l$ of
the generators $a_i$.

\subsection{The geometric model of $\G(S)$}

When $\G = \G\s$ and the $n_i$ are relatively prime, but not all
prime, the product of trees on which $\G = \G\s$ acts is not as
simple.

We first discuss the tree $T^n$ on which the group $BS(1,n)$ acts
when $n$ is not prime.  Suppose that $n = \pexp$, and let $T_i$ be
the Bruhat-Tits tree associated to $PSL_2(\Q_{p_i})$.  The tree
$T^n$ is a subspace of $\ptr$, whose branching may not be constant,
but depends on the exponents of the primes.

Define a folding function $F_i:T_i \rightarrow \R$ as follows.  If
$h_i$ is the height function defined on $T_i$, then $F_i(t) =
h_i(t)$ for $t \in T_i$.  Combining folding functions on the
$T_i$ we get a map $F(r): \ptr \rightarrow \R^r$ defined by $F(r)
= (F_1,F_2, \cdots ,F_r)$.  Consider the grid of lines in $\R^r$
of the form $(x_1,x_2, \cdots, x_{j-1}, \R , x_{j+1}, \cdots
,x_r)$ where $x_i \in \Z$.  So we really have $r$ families of
parallel lines in $\R^r$.  View each family as representing folded
copies of one of the trees $T_i$ under $F(r)$.

The branching of the tree $T^n$ is determined by the line $e_1x_1
+ e_2x_2 + \cdots + e_r x_r  = 0$ in $\R^r$.  When the line
crosses a line in the family of parallel lines corresponding to
$T_i$, the tree $T^n$ branches $n$ times.  When the line crosses
the intersection of two lines, one from the family of $T_i$ and
one from the family of $T_j$, the branching is $i+j$.

\medskip
\noindent {\bf Example.} Consider the group $BS(1,6)$.  The tree
$T^6$ on which it acts is a subset of $T_2 \times T_3$, determined
by the line $y=x$ in the plane $\R^2$, since the exponent of each
prime is $1$.  This line only crosses vertices of the grid of
lines, so the branching is uniform of valence $6$.

\medskip
\noindent {\bf Example.} Consider the group $BS(1,12)$.  The tree
$T^{12}$ on which it acts is a subset of $T_2 \times T_3$, only
now the line in $\R^2$ which determines the branching is $2x=y$.
From the way this line crosses the grid of lines, we see that the
branching of $T$ is not uniform.  The vertices alternate between
valence $2$ and valence $6$, where the valence $2$ arises from the
line crossing only a horizontal grid line, and the valence $6$
arises when the line crosses a vertex in this grid of lines.

\medskip
\noindent {\bf Example.} Consider the group $BS(1,60)$.  The tree
$T^{60}$ on which it acts is a subset of $T_2 \times T_3 \times
T_5$, and now the folding map $F(3)$ is a map to $\R^3$.  The line
in $\R^3$ determining the branching of $T$ is $2x+y+z = 0$.  Again
we see that the amount of branching at each valence varies.

\medskip
Now consider $\G = \G\s$ where the $n_i$ are not all prime.
Consider any $n_i$, and let $\p$ be the list of primes dividing
$n_i$, with $T_i$ the Bruhat-Tits tree of $PSL_2(\Q_{p_i})$.
Then let $T^i$ be the tree on which $BS(1,n_i)$ acts
(described above) which is a subspace of $\ptr$.  Then $\G$ acts
on $\pti$.   The complex $X(S)$ is then warped product of $\pti$ with $\R$.

\section{The Structure of Quasi-Isometries}
\label{sec:proofs}

The key step in the proofs of the theorems in this paper is understanding the
structure of the quasi-isometries of $\G$, or equivalently of $X$.
We begin with two groups $\G_1$ and $\G_1$ and a $(K,C)$-quasi-isometry between
their geometric models, $f: X_1 \rightarrow X_2$.

Let $\pi$ be the projection $X \to \pt$.  Define a {\em horocycle} of
the complex $X$ to be a subset of the form $\pi^{-1}(\point)$ where $\point$ is a point in $\pt$.  A {\em hyperplane} in $X$ is a subcomplex of the form $\pi^{-1}( \times l_1 \times \cdots \times l_n)$
where each $l_i$ is a geodesic in $T_i$. The first goal is to show that the
quasi-isometry $f$ preserves horocycles and hence induces a quasi-isometry of a
product of trees.  These arguments are similar to those in \cite{W2}.

\begin{lemma}
\label{lemma:SW}
For any $(K,C)$ there is an $R>0$ so that for any $f: X_1 \to X_2$, a
$(K,C)$-quasi-isometry, and every hyperplane $H$ of $X_1$, there is a subset $Y$ of
$\pt$ so that the image $f(H)$ is within $R$ of $\pi^{-1}(Y)$.
\end{lemma}

\begin{proof}

Let $g$ be a quasi-inverse of $f$, so that $g \circ f$ is a bounded
distance from the identity map, and hence proper. By a standard connect-the-dots
argument, we may assume $f$ and $g$ are continuous.  As the $X_i$ are uniformly
contractible, the compositions are homotopic to the identity through homotopies
of length at most $R_0$ (depending only on the $X_i$ and the constants $(K,C)$).   The the
maps $f$ and $g$ are, in particular, proper homotopy equivalences.

Consider the fundamental class $[H]$ in $H^{uf}_{n+1}(X_1)$.  The push forward
$f_*([H])$ is thus a non-trivial class in $H^{uf}_{n+1}(X_2)$.  Further, this
class clearly has a representative $c$ with support contained in the $R_0$
neighborhood of $f(H)$.  The simplicial structure of $X_2$ forces the coefficients
of $c$ to be constant along horocycles.  Thus the support of $c$ is of the form
$Y \times \R$ for some subcomplex $Y$ of $\pt$.   This shows that the $R_0$-neighborhood of
$f(H)$ contains $Y \times \R$.

To complete the proof we must show that a neighborhood of $Y \times \R$ ($= \hbox{supp}(c)$)
contains $f(H)$.  If not, then there are arbitrarily large balls in $f(H)$ which are
not contained in $Y \times \R$.  Applying the inverse map, $g$, this would give a
representative of $[H]$ whose support misses large balls in $H$.  This is impossible,
as any representative of the fundamental class has full support.
\end{proof}

\begin{lemma}[Horocycles are preserved]
\label{lemma:horocycles}
For every $(K,C)$ there is an $R$ so that if $f$ is a $(K,C)$-quasi-isometry of $X$ and $h$ is a horocycle in $X$ then there is a horocycle $h'$ so that $d_H(f(h),h')\leq R$.
\end{lemma}

\begin{proof}

This is an immediate consequence of the previous lemma.  For any horocycle $h$
there are a finite number of hyperplanes $H_1, \cdots, H_k$ in $X_1$ which have
coarse intersection at Hausdorff distance at most $R$ from $h$ (the constants $k$ and $R$ depend only on the geometry of $X_1$).  The previous lemma implies that the image of $h$
is Hausdorff equivalent to a complex of the form $Y \times \R$ for some subset $Y$
of $X_2$.  Applying the same argument to the inverse map $g$ and each horocycle
in $Y \times \R$, we conclude that $Y$ must be of finite diameter (bounded
independently of $h$).  This proves the lemma.
\end{proof}

\begin{corollary}[Factor preserving]
\label{cor:samenumber} Consider the groups $\Gamma_1 = \Gamma\s$
and $\Gamma_2 = \Gamma\sprime$ where $(n_i,n_j) = (m_i,m_j) = 1$ for $i \neq j$, and a quasi-isometry $f: \G_1
\rightarrow \G_2$ between them.  Then:

\begin{enumerate}
\item
$k=l$,

\item
$f$ induces a quasi-isometry $f_T: \pt \rightarrow \ptp$ ,and

\item
there are $(K',C')$-quasi-isometries $\overline{f_i}: T_i \rightarrow T'_i$ (after possibly reindexing the
tree factors) so that $f_T$ is a bounded distance from the product quasi-isometry
$\overline{f_1} \times \cdots \times \overline{f_k}$.
\end{enumerate}
\end{corollary}

\begin{proof}
Since every point $(t_1,t_2, \cdots ,t_n) \in \pt$ determines a
horocycle in $X$, it follows from lemma \ref{lemma:horocycles}
that the quasi-isometry $f$ induces a quasi-isometry on the
product of trees: $f_T:\pt \rightarrow \ptl$. It now follows from
theorem \ref{thm:KL} that $k=l$ and thus there are the
same numbers of parameters in $\G_1$ and $\G_2$.  It then follows
from theorem \ref{thm:KL}  that this map is a bounded
distance from a product $f_1 \times \cdots \times f_k$ of
quasi-isometries.
\end{proof}

\begin{corollary}[Bilipschitz maps]
\label{cor:factorpreserving}
Let $f: \Gamma_1 \rightarrow \Gamma_2$ be a $(K,C)$-quasi-isometry.
Then there are bilipschitz maps $g: \R \rightarrow \R$ and $fi: T_i
\rightarrow T_i$ (after possibly re-indexing the tree factors)
so that $f$ is a bounded distance from $(g,f_1, \cdots f_k)$.
\end{corollary}

\begin{proof}
Applying corollary \ref{cor:samenumber} we may assume that the quasi-isometry $f$ preserves the individual tree factors.  We use the notation of corollary \ref{cor:samenumber} and let $f_i$ denote the induced map on the $i$-th tree factor.  It follows that the
quasi-isometry $f$ restricts to a map on each Baumslag-Solitar
subcomplex $T_i \times \R$, which is also a quasi-isometry.
Applying theorem \ref{thm:BS:2} of Farb and Mosher, we conclude that $f_i$ is a bounded distance from the product of a bilipschitz map of $T_i$ with a bilipschitz map of $\R$.  It is easy to see that we must obtain the same bilipschitz map of $\R$ regardless of which Baumslag-Solitar subspace we restrict to, and the corollary follows.
\end{proof}

We are now able to prove theorem \ref{thm:classification}.

\medskip
\noindent {\it Proof of Theorem \ref{thm:classification}.}
Applying corollary \ref{cor:factorpreserving}, we consider our
quasi-isometry to be factor preserving of the form $(g,f_1, \cdots
f_n)$, with each individual map bilipschitz. Then any pair
$(g,f_i): \R \times T_i \rightarrow \R \times T_i'$ is a
quasi-isometry of $BS(1,p_i)$ to $BS(1,q_i)$, by theorem
\ref{thm:BS:2} . It follows from theorem \ref{thm:BS:1}
that, after reordering, $p_i = q_i$.
\qed

\subsection{Description of the quasi-isometry group}
\label{sec:QIgroup}

We begin with a lemma important for the proof of theorem \ref{thm:QIgroup}.

\begin{lemma}[\cite{FM2}Rubber Band Principle]
\label{lemma:rubberband}
For all $L,M >0$ there is a constant $C$ satisfying the following
property.
Suppose $X$ and $Y$ are path metric spaces and $f:X \rightarrow Y$ is
a map.  Suppose that there are collections of isometrically embedded
subspaces $C_X$ of $X$ and $C_Y$ of $Y$ satisfying:

\begin{itemize}
\item
Any two points in $X$ (or in $Y$) can be connected by an
$M$-quasi-geodesic made up of a finite number of subpaths, each lying
in an element of $C_X$ or $C_Y$.

\item
$f$ induces a one-to-one correspondence between elements of $C_X$ and
$C_Y$.

\item
$f$ restricts to an $L$-quasi-isometry between corresponding elements
of $C_X$ and $C_Y$.

\end{itemize}
Then $f:X \rightarrow Y$ is a $C$-quasi-isometry.

\end{lemma}

We are now able to prove theorem \ref{thm:QIgroup} and describe the
quasi-isometry group of $\Gamma$.

\bigskip
\noindent
{\it Proof of Theorem \ref{thm:QIgroup}.}
It is clear that we have a homomorphism $$\Phi: QI(X) \rightarrow Bilip(\R) \times Bilip(\Q_{p_1}) \times Bilip(\Q_{p_1}) \times \cdots
\times Bilip(\Q_{p_k})$$ given by
$
\Phi(f) = (f_R,f_1,f_2, \cdots ,f_n).$
In addition we get a homomorphism $$\Phi_i: QI(X) \rightarrow
QI(BS(1,p_i)) \cong Bilip(\R) \times
Bilip(\Q_{p_i})$$
for each $i = 1,2, \cdots ,n$, given by $\Phi(f) = (f_R,f_i)$.
Following the reasoning in \cite{FM}, we see that for any $f \in
ker(\Phi)$, the quasi-isometry $\Phi_i(f)$ is a bounded distance $B_i$
from the identity map on $X_n$.  Letting $B = \max \{B_i \}$, the
Rubber Band Principle implies that $\Phi(f)$ is a bounded distance $B$
from the identity.

To see that $\Phi$ is surjective, we again use the Rubber band
Principle to piece together quasi-isometries of the $X_n$
subcomplexes.
Choose $f_R \in Bilip(\R)$ and maps $f_i \in Bilip(\Q_{p_i})$.  We must
show that $f_R \times f_1 \times \cdots \times f_n$ is a
quasi-isometry of $\Gamma_1$.
From \cite{FM} we know that any pair $(f_R,f_i)$ yields a
quasi-isometry of $X_i$.
We can assume the quasi-isometry constants are uniform by taking the
largest pair of constants from any of these maps.
Since the $f_R$ is common to any two pairs,  we obtain a product map
$f$ of the
entire complex.
Thus we have a collection of subspaces and map $f = f_R \times f_1 \times
\cdots \times f_n$ satisfying the Rubber Band Principle, so $f$
is a quasi-isometry of $\Gamma$.
\qed

\section{Rigidity}

We finish with the proof of theorem \ref{thm:rigidity}, which shows that this class
of groups is quasi-isometrically rigid.

\medskip
\noindent
{\it Proof of Theorem \ref{thm:rigidity}.}
Let $\G'$ be any finitely generated group quasi-isometric to $\G(S)$, with $X$ a the
model space for $\G(S)$ as before, and let $f: \G' \to X$ be a quasi-isometry with $g$
a coarse inverse.  We get a quasi-action of $\G'$ on $X$ where $\g' x$=$f(\g' g(x))$.
By lemma \ref{lemma:horocycles}, horocycles are preserved, so we get an induced quasi-action of $\G'$ on the
product of trees $\pt$.  By passing to a finite index subgroup of $\G'$ we may assume that
the quasi-action is the diagonal quasi-action of a collection of quasi-actions $\G'$ on
$T_i$.  The maps to the complexes of the Baumslag-Solitar
subgroups of $\G(S)$ are $\G'$ equivariant (to within finite distance), and so quasi-preserve the height function.  By \cite{MSW}, there are trees $T_i'$ quasi-isometric to the $T_i$ and actions of $\G'$ on $T_i$ which are quasi-conjugate to the quasi-actions on the $T_i$.  Further, each of these trees is  homogeneous with a $\G'$ invariant orientation with one edge directed into each vertex.
Thus we get an action of $\G'$ on the product of the $T_i'$, with vertex stabilizers
virtually cyclic, preserving the orientations, and with finite quotient, in other words
we get a description of $\G'$ as a finite complex of virtual $\Z$s.

 Consider the edges in this quotient which come from edges of a $T_i'$. These are oriented,
and as there is exactly one edge oriented toward every vertex in $T_i'$, the same is true
in the quotient.  Since the quotient is finite, this implies that there is precisely one
such edge oriented away from each vertex of the quotient.  This implies that these edges
consist of a finite union of circles.  Further, this implies that for any $v \in T_i'$, the
action of $stab(v)$ on edges pointing away from $v$ is transitive. Similarly, fixing any
edge in $T_i'$ and looking at two cells coming from $T_i' \times T_j'$ we have exactly two
such two cells at every edge of the quotient, with one oriented towards, and one away from,
this edge.  Continuing over higher dimensional cubes, we see that the quotient is product of
oriented circles, with the inclusions of the stabilizer of a cube to a face stabilizer is
an isomorphism if it goes against the orientation.  Thus we may collapse such a cube, unless
its opposite faces are the same in the quotient.   Making all such possible collapses gives
a complex of groups description of $\G'$ with underlying complex a product of oriented loops,
with stabilizers all virtually $\Z$, and with the inclusions isomorphism when going against
the orientation.  As in \cite{FM}, we may pass to a finite index subgroup of $\G'$ which has
such a description with all stabilizers $\Z$.  Such a complex of groups has a presentation
precisely of the form $\G(S')$ for some $S'$.   Thus $\G'$ is commensurable to $G(S')$ for
some $S'$, as desired.
\qed

\medskip
\noindent
Jennifer Taback\\
Dept. of Mathematics and Statistics\\
University of Albany\\
Albany, NY 12222\\
jtaback@math.albany.edu

\bigskip
\noindent
Kevin Whyte\\
Dept. of Mathematics\\
University of Illinois at Chicago\\
Chicago, Il 60607\\
kwhyte@math.uchicago.edu

\end{document}